\documentclass{article}

\usepackage{amssymb}
\usepackage{amsmath}

\usepackage[dvips]{graphicx}
\begin{document}

%%% remove comment delimiter ('%') and select language if required
%\selectlanguage{spanish} 

\begin{center}
\noindent \textbf{A non-unital algebra has UUNP iff its unitization has UUNP}                                                                                                                                                                                     
\end{center}

\noindent \textbf{}

\begin{center}
M. El Azhari
\end{center}

\noindent \textbf{ } 

\noindent \textbf{Abstract.} Let $A$ be a non-unital Banach algebra, S. J. Bhatt and H. V. Dedania showed that $A$ has the unique uniform norm property (UUNP) if and only if its unitization has UUNP. Here we prove this result for any non-unital algebra.

\noindent \textbf{}


\noindent \textbf{}


\noindent \textbf{} 
 
\noindent \textbf{} Let $A$ be a non-unital algebra and let $A_{e}=\lbrace a+\lambda e: a\in A,\,\lambda\in C\rbrace$ be the unitization of $A$ with the identity denoted by $e.$ For an algebra norm $\Vert .\Vert$ on $A,$ define $\Vert a+\lambda e\Vert_{op}=\sup \lbrace\Vert(a+\lambda e)b\Vert: b\in A, \Vert b\Vert\leq 1\rbrace $ and $\Vert a+\lambda e\Vert_{1}=\Vert a\Vert +\vert\lambda\vert$ for all $a+\lambda e\in A_{e}.\:\Vert .\Vert_{op}$ is an algebra seminorm on $A_{e},$ and $\Vert .\Vert_{1}$ is an algebra norm on $A_{e}.$ An algebra norm $\Vert .\Vert$ on $A$ is called regular if $\Vert .\Vert_{op}=\Vert .\Vert$ on $A.$ A uniform norm $\Vert .\Vert$ on $A$ is an algebra norm satisfying the square property $\Vert a^{2}\Vert = \Vert a\Vert^{2}$ for all $a\in A;  $ and in this case, $\Vert .\Vert$ is regular and $\Vert .\Vert_{op}$ is a uniform norm on $A_{e}.$ An algebra has the unique uniform norm property (UUNP) if it admits exactly one uniform norm.

\noindent \textbf{}

\noindent \textbf{Theorem .} A non-unital algebra $A$ has UUNP if and only if its unitization $A_{e}$ has UUNP.

\noindent \textbf{}

\noindent \textbf{}Proof: Let $\Vert .\Vert$ and $ \vert\vert\vert .\vert\vert\vert$ be two uniform norms on $A_{e},$ then $\Vert .\Vert =\vert\vert\vert .\vert\vert\vert$ on $A$ since $A$ has UUNP, and so $\Vert .\Vert_{op} =\vert\vert\vert .\vert\vert\vert_{op}$ on $A_{e}.$ By [3, Corollary 2.2(1)] and since two equivalent uniform norms are identical, it follows that $(\Vert .\Vert =\Vert .\Vert_{op}$ or $\Vert .\Vert \cong\Vert .\Vert_{1})$ and $( \vert\vert\vert .\vert\vert\vert =\vert\vert\vert .\vert\vert\vert_{op} =\Vert .\Vert_{op}$ or $\vert\vert\vert .\vert\vert\vert \cong\vert\vert\vert .\vert\vert\vert_{1} =\Vert .\Vert_{1});$ equivalently, at least one of the following holds:\\
(i) $\Vert .\Vert =\Vert .\Vert_{op}$ and $ \vert\vert\vert .\vert\vert\vert =\vert\vert\vert .\vert\vert\vert_{op} =\Vert .\Vert_{op};$ \\
(ii) $\Vert .\Vert =\Vert .\Vert_{op}$ and $\vert\vert\vert .\vert\vert\vert \cong\vert\vert\vert .\vert\vert\vert_{1} =\Vert .\Vert_{1};$ \\
(iii) $\Vert .\Vert \cong\Vert .\Vert_{1}$ and $ \vert\vert\vert .\vert\vert\vert =\vert\vert\vert .\vert\vert\vert_{op} =\Vert .\Vert_{op};$ \\
(iv) $\Vert .\Vert \cong\Vert .\Vert_{1}$ and $\vert\vert\vert .\vert\vert\vert \cong\vert\vert\vert .\vert\vert\vert_{1} =\Vert .\Vert_{1}.$ \\
If either (i) or (iv) is satisfied, then $\Vert .\Vert =\vert\vert\vert .\vert\vert\vert.$ By noting that (ii) and (iii) are similar by interchanging the roles of $\Vert .\Vert$ and $\vert\vert\vert .\vert\vert\vert,$ it is enough to assume (ii). Let $(c(A),\Vert .\Vert^{\sim})$ be the completion of $(A,\Vert .\Vert),$ we distinguish two cases:\\
(1) $c(A)$ has not an identity:\\
$\Vert .\Vert^{\sim}$ is regular since it is uniform. By [1, Corollary 2], $\Vert .\Vert^{\sim}_{op}\leq\Vert .\Vert^{\sim}_{1}\leq 3\Vert .\Vert^{\sim}_{op}$ on $c(A)_{e}$ (unitization of $c(A)$). Let $a+\lambda e\in A_{e}\subset c(A)_{e},\: \Vert a+\lambda e\Vert^{\sim}_{1}=\Vert a\Vert^{\sim}+\vert\lambda\vert=\Vert a\Vert+\vert\lambda\vert=\Vert a+\lambda e\Vert_{1}$ and $\Vert a+\lambda e\Vert^{\sim}_{op}=\sup\lbrace\Vert (a+\lambda e)b\Vert^{\sim}: b\in c(A), \Vert b\Vert^{\sim}\leq 1\rbrace =\sup\lbrace\Vert (a+\lambda e)b\Vert: b\in A, \Vert b\Vert\leq 1\rbrace =\Vert a+\lambda e\Vert_{op}.$ Therefore $\Vert .\Vert_{op}\leq\Vert .\Vert_{1}\leq 3\Vert .\Vert_{op}.$ By (ii), $\Vert .\Vert$ and $\vert\vert\vert .\vert\vert\vert$ are equivalent uniform norms, and so $\Vert .\Vert =\vert\vert\vert .\vert\vert\vert.$\\
(2) $c(A)$ has an identity $e:$\\
Let $(c(A_{e}),\vert\vert\vert .\vert\vert\vert^{\sim})$ be the completion of $(A_{e},\vert\vert\vert .\vert\vert\vert).$ Since $\Vert .\Vert =\vert\vert\vert .\vert\vert\vert$ on $A,\:c(A)$ can be identified to the closure of $A$ in $(c(A_{e}),\vert\vert\vert .\vert\vert\vert^{\sim})$ so that $\Vert .\Vert^{\sim} =\vert\vert\vert .\vert\vert\vert^{\sim}$ on $c(A).$ Let $a+\lambda e\in A_{e}\subset c(A),$\\
$\Vert a+\lambda e\Vert =\Vert a+\lambda e\Vert_{op}$ by (ii)\\
$=\sup\lbrace\Vert (a+\lambda e)b\Vert: b\in A, \Vert b\Vert\leq 1\rbrace\\
=\sup\lbrace\Vert (a+\lambda e)b\Vert^{\sim}: b\in c(A), \Vert b\Vert^{\sim}\leq 1\rbrace \\
= \Vert a+\lambda e\Vert^{\sim}$ since $c(A)$ is unital\\
$=\vert\vert\vert a+\lambda e\vert\vert\vert^{\sim}=\vert\vert\vert a+\lambda e\vert\vert\vert.$ Thus $\Vert .\Vert =\vert\vert\vert .\vert\vert\vert.$\\
Conversely, let $\Vert .\Vert$ and $\vert\vert\vert .\vert\vert\vert$ be two uniform norms on $A,$ then $\Vert .\Vert_{op}$ and $\vert\vert\vert .\vert\vert\vert_{op}$ are uniform norms on $A_{e},$ hence $\Vert .\Vert_{op}=\vert\vert\vert .\vert\vert\vert_{op}$ since $A_{e}$ has UUNP. Therefore $\Vert .\Vert =\Vert .\Vert_{op}=\vert\vert\vert .\vert\vert\vert_{op}=\vert\vert\vert .\vert\vert\vert$ on $A$ since $\Vert .\Vert$ and $\vert\vert\vert .\vert\vert\vert$ are regular.

\noindent \textbf{}

\noindent \textbf{References}

\noindent \textbf{}

\noindent \textbf{}[1] J. Arhippainen and V. M\"{u}ller, Norms on unitizations of Banach algebras revisited, Acta Math. Hungar., 114(3)(2007), 201-204.

\noindent \textbf{}[2] S. J. Bhatt and H. V. Dedania, Uniqueness of the uniform norm and adjoining identity in Banach algebras, Proc. Indian Acad. Sci. (Math. Sci.), 105(4)(1995), 405-409.

\noindent \textbf{}[3] H. V. Dedania and H. J. Kanani, A non-unital $\ast$-algebra has $UC^{\ast}NP$ if and only if its unitization has $UC^{\ast}NP$, Proc. Amer. Math. Soc., 141(11)(2013), 3905-3909.
 
\noindent \textbf{}
 
\noindent \textbf{}

\noindent \textbf{}Ecole Normale Sup\'{e}rieure, Avenue Oued Akreuch,   Takaddoum, BP 5118, Rabat,  Morocco.
  
\noindent \textbf{}

\noindent \textbf{}E-mail: mohammed.elazhari@yahoo.fr

\end{document}